\documentclass[12pt,twoside,a4paper,notitlepage]{article}

\usepackage{amsmath}

\usepackage{amssymb}
\usepackage{amsfonts}
\usepackage[a4paper,margin=2cm,vmargin={1cm,2cm},includeheadfoot]{geometry}
\usepackage[T1]{fontenc}
\usepackage{amsthm}
\usepackage{enumerate}
\usepackage[scaled=0.92]{helvet}
\usepackage{srcltx}
\usepackage{graphicx}
\usepackage{tikz}
\usepackage{booktabs}

\usepackage{fancybox}
\usepackage{textcomp}

\usepackage{fancyvrb}

\usepackage[labelfont=bf,margin=1cm,justification=justified,labelsep=period]{caption}

\theoremstyle{plain}
\newtheorem{theorem}{Theorem}
\newtheorem{proposition}[theorem]{Proposition}

\newtheorem{lemma}[theorem]{Lemma}

{\theoremstyle{definition}
\newtheorem{remark}[theorem]{Remark}}

{\theoremstyle{definition}
}

\newcommand{\R}{{\mathbb {R}}}

\newcommand{\C}{{\mathbb{C}}}

\newcommand{\E}{\text{\sf \bfseries E}}
\newcommand{\p}{\text{\sf \bfseries P}}
\newcommand{\Prob}{\text{\sf \bfseries P}}

\newcommand{\bs}{\boldsymbol}   
\newcommand{\wt}{\widetilde}    
\newcommand{\wh}{\widehat}      
\newcommand{\1}{{\bf 1}}        

\newcommand{\buil}[3]{\mathrel{\mathop{\kern0pt#1}\limits_{#2}^{#3}}}

\begin{document}

\begin{center}
{\bf \large Efficient computation of first passage times in Kou's
jump-diffusion model} \\ [8 pt]

{\sc Abdel Belkaid}\footnote{{\it Department of Production,
 Technology and Operations Management, IESE Business School, Barcelona Campus,
 Avda. Pearson, 21, 08034 Barcelona, Spain},
  e-mail: {\tt abelkaid@iese.edu}}
and
{\sc Frederic Utzet}\footnote
{Corresponding author. {\it Department of Mathematics, Universitat Autònoma de Barcelona,
Campus de Bellaterra,
08193 Bellaterra (Barcelona), Spain},
e-mail: {\tt utzet@mat.uab.cat}}


\end{center}

\begin{center}
\begin{minipage}{0.8\linewidth}
\begin{small}

\noindent{\bf Abstract.}
S. G.  Kou and H. Wang [First Passage times of a Jump Diffusion Process
  \textit{Ann.  Appl. Probab.} {\bf 35} (2003) 504--531]  give expressions of both the (real)
Laplace transform of the distribution of  first passage time
and the (real) Laplace transform of the joint distribution of the first passage
 time and the running maxima of a jump-diffusion model  called Kou model.
   Kuo and Wang invert the first
  Laplace transform by using Gaver-Stehfest algorithm, and  the inversion of
  second one involves a large computing time with an algebra computer system.
 In the present paper,  we
give  a much  simpler expression of the Laplace transform of the joint
 distribution, and  we also show, using Complex Analysis techniques, that
both Laplace transform can be extended to the complex plane. Hence,
we can use variants of the Fourier-series methods to invert that Laplace transfoms,
which are very efficent.
 The improvement in the computing times and  accuracy is remarkable.

\smallskip


\smallskip


\end{small}
\end{minipage}
\end{center}

\section{Introduction}
The Kou model (Kou 2000, and Kou and Wang 2003)
 is a jump-diffusion model with the jumps given by a  mixture of
a positive and a negative exponential random variables.
It is one of the few L\'evy processes, allowing  positive and negative
jumps, which permits the probability of first passage times to be computed
analytically.
This property, jointly with the fact that it is very flexible when
 applied to real data,  makes the Kou model very interesting in applications.
 In   Mathematical Finance, in particular, the pricing of important
 derivative products, as barrier or lookback options,  relies on
 these  computations. We should point out that the Kou model belongs to
 the family of meromorphic L\'evy processes (see
Kuznetsov 2012  and Kuznetsov {\it et al.} 2011).

Kou and Wang (2003) obtained a closed expression of the real Laplace
 transform
of the distribution of the  first passage time, and an expression of the
real Laplace transform of the joint
distribution
of the process and the running maxima in terms of the so-called {\sf Hh} functions.
 They also  proposed ways to invert
these Laplace transforms. However, for the inversion of the Laplace
transform of the distribution of the first passage
time they use Gaver-Stehfest algorithm, which has the advantage that
it only uses
real numbers but, in general, it is difficult to control its accuracy.
For  the joint distribution, the presence of the {\sf Hh} functions  means that
 the inversion procedure  needs a large  computation time
with a software package like MATHEMATICA. 
In the present paper, we first show that the Laplace transform of the first
passage time can be extended to the complex plane and can be inverted by using the
 complex inversion formula or  Bromwich integral. Second,  we give a closed form of the
  Laplace transform
of the joint distribution  of the process and the running maxima without {\sf Hh}
 functions, and
  also prove that it can be inverted by Bromwich integral. As a consequence,
  very efficient methods can be applied to invert both Laplace transforms.
  As regards  the times
reported by Kou and Wang (2003), the computation time
of the joint distribution is reduced from minutes to milliseconds.

In Appendix 2 there are the codes of some of the Maple and C programs used by the authors.

\section{Kou model: notations}
\label{sec_Kou}
\begin{center}
\begin{figure}[htb]
\begin{center}
\includegraphics[width=10cm]{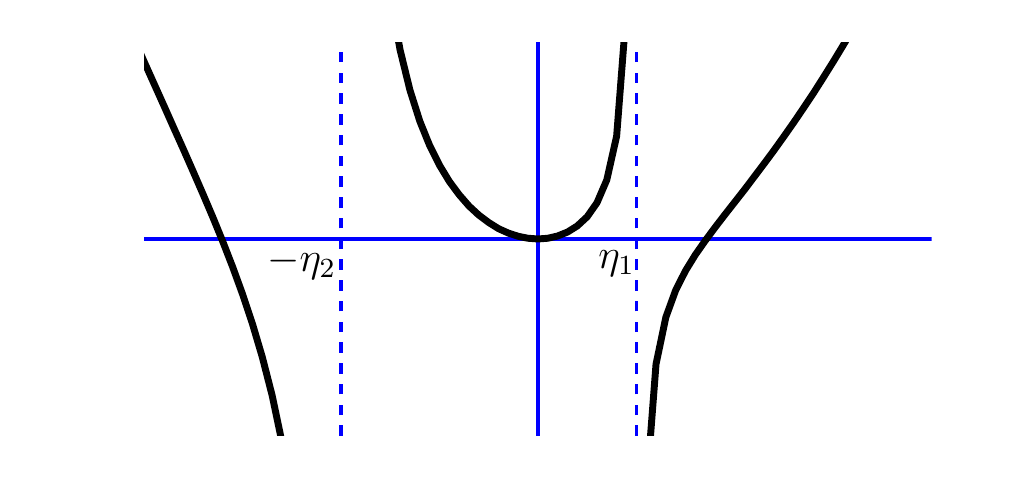}
\end{center}
\caption{Typical graph of  the function $G(x)$}
\label{kuo:fig}
\end{figure}
\end{center}

We use the same notations as Kou and Wang (2003).
 The Kou model   is  a  L\'evy process of the type jump-diffusion of the form
 \begin{equation}
\label{hyper}
X_t=\sigma W_t+\mu t+\sum_{j=1}^{N_t} Y_j,
\end{equation}
where $\mu\in \R$,  $W=\{W_t,\ t\ge 0\}$ is a standard Brownian motion,
 $\sigma>  0 $,
$N=\{N_t,\ t\ge 0\}$ is a Poisson process of intensity $\lambda>0$, independent
of $W$, and $\{Y_n,\ n\ge 1\}$ is a sequence of i.i.d. random variables,
 independent
of $W$ and $N$, and with a law given by a  mixture of a positive and a negative
 exponential random variables of parameters $\eta_1>0$ and $\eta_2>0$ respectively;
   that is, the density of $Y_n$ is
 $$f_Y(y)=p \eta_1 e^{-\eta_1 y} \1_{(0,\infty)}(y)+
(1-p) \eta_2 e^{\eta_2 y} \1_{(-\infty,0)}(y),$$
where $0<p<1$.
The characteristic function of $X_t$ is
$$\E\big[e^{isX_t}\big]=e^{tG(is)},\ s\in \R,$$
where
\begin{equation}
\label{G}
G(x)=\mu x+\frac{1}{2}\sigma^2 x^2+\lambda\Big(\frac{p\eta_1}{\eta_1-x}+
\frac{(1-p)\eta_2}{\eta_2+x}-1\Big).
\end{equation}
Kou and Wang (2003) Lemma 2.1 prove that the for every $\alpha>0$,
the equation $G(x)=\alpha$ has exactly four real zeros: $\beta_{1,\alpha},
\beta_{2,\alpha},-\beta_{3,\alpha}$ and $-\beta_{4,\alpha}$,
where
$$0<\beta_{1,\alpha}<\eta_1<\beta_{2,\alpha}
\quad \text{and}\quad 0<\beta_{3,\alpha}<\eta_2<\beta_{4,\alpha}.$$
See in Figure \ref{kuo:fig} a typical plot of the function $G(x)$.

We  also have the  following equality:
\begin{equation}
\label{fact_q-psi}
\frac{\alpha}{\alpha-G(x)}=\frac{1-x/\eta_1}{(1-x/\beta_{1,\alpha})(1-x/
\beta_{2,\alpha})}\,
\frac{1+x/\eta_2}{(1+x/\beta_{3,\alpha} )(1+x/\beta_{4,\alpha})}.
\end{equation}
This equality  also holds for $x=is,\ s\in \R$, and indeed
 it is the Wiener-Hopft factorization related to the L\'evy process (see Sato 1999).

\section{Theoretical results}
\subsection{First passage times}
 For a level $b> 0$
  denote by $\tau_b$ the first  passage time above $b$:
 $$\tau_b=\inf\{t\ge 0:\, X_t \ge b\}.$$
 We have that $\tau_b\le t$ if and only if $\max_{0\le s\le t} X_s\ge b.$
Our purpose is to compute  the probability $\p(\tau_b\le t)$.
To this end,
 Kou and Wang (2003), Theorem 3.1 and formula on line 16 of page 512,
 obtained the Laplace transform
  of that  probability (see also Lemma \ref{lema:tranformada:p} in the Appendix):

\begin{proposition}(Kou and Wang  2003)
\label{Prop:teoKou} Fix $b>0$.
 For every $\alpha>0$,
\begin{equation}
\label{laplace:uni}
\int_0^\infty e^{-\alpha t}\,\p(\tau_b\le t)\, dt=
\frac{1}{\alpha}\bigg(\frac{\beta_{2,\alpha}\big(\eta_1-\beta_{1,\alpha})}{\eta_1\big(
\beta_{2,\alpha}-\beta_{1,\alpha}\big)}\, e^{-b \beta_{1,\alpha}}
+\frac{\beta_{1,\alpha}\big(\beta_{2,\alpha}-\eta_1\big)}
{\eta_1\big(\beta_{2,\alpha}-\beta_{1,\alpha}\big)}\, e^{-b \beta_{2,\alpha}}\bigg),
\end{equation}
where $0<\beta_{1,\alpha}<\beta_{2,\alpha}$ are the positive zeros of the equation
$G(x)=\alpha.$
\end{proposition}

We prove in the Appendix that this equality can be extended to the
 half complex plane $\{\alpha\in \C:\, {\rm Re}(\alpha)>0\}$ except on a finite number of removable singularities.
The key point is that by analytic continuation
we can construct four analytic functions (of $\alpha\in \C$), that we denote
as before by $\beta_{1,\alpha},
\beta_{2,\alpha}, -\beta_{3,\alpha},-\beta_{4,\alpha}$ which
are the  four (complex) roots of the equation $G(z)=\alpha$, for
$\alpha\in\C$.
We prove that for $\alpha\in \C$, with $\text{Re}(\alpha)>0$,
two of these functions always have  positive real part and two have negative real
part, and we provide a procedure to compute the singular points.

More specifically, $G(z)$ given in (\ref{G}) is a rational function with numerator of degree 4 and denominator of degree
2, and we can write (see (\ref{fact_q-psi}))
 \begin{equation}
\label{GP}
G(z)-\alpha=\frac{P_\alpha(z)}{Q(z)}
\end{equation}
where $Q(z)=(\eta_1-z)(\eta_2+z)$, and $P_\alpha(z)$ is a fourth degree polynomial
whose coefficients are polynomials on $\alpha$ (see an example in Section \ref{sec_Num}), and
its roots are the zeros of the equation $G(z)=\alpha$.

In principle,  the function defined by the right hand side of (\ref{laplace:uni})
has  singular points (besides   $\alpha=0$)  at the values of
 $\alpha$ where $\beta_{1,\alpha}=\beta_{2,\alpha},$
that is, where
$P_\alpha(z)$ has a multiple root. However,  these points
are the roots of
 the resultant $R(\alpha)$ of $P_\alpha(z)$ and $\partial P_\alpha/\partial z$,
 which is a polynomial in $\alpha$.
 The behaviour of the function at these points  can be controlled.
 See the proof of Proposition \ref{prop_lap_comp1}  for details and references.

\begin{proposition}
\label{prop_lap_comp1} For  $\alpha \in \C$ with $\rm{Re}(\alpha)>0$, the equality (\ref{laplace:uni})
holds,
except at some
removable  singularities given by the roots of the resultant of the polynomials
$P_\alpha(z)$ and $\partial P_\alpha/\partial z$, where $P_\alpha(z)$ is given by
(\ref{GP}),  and $\beta_{1,\alpha}$ and $\beta_{2,\alpha}$ are
 the zeros of $G(z)=\alpha$ with ${\rm Re}(\beta_{1,\alpha})>0$
 and ${\rm Re}(\beta_{2,\alpha})>0$.

\end{proposition}

We prove that the complex Laplace transform (\ref{laplace:uni}) can be inverted by using Bromwich integral. Specifically,

\begin{proposition}
\label{laplace:inv_complexa1} For any  $b>0$ and $t>0$,
\begin{equation}
\label{laplace:inversio}
\p(\tau_b\le t)=\frac{1}{2\pi i}\int_{u-i\infty}^{u+i\infty}
\frac{1}{\alpha}\!\bigg(\!\frac{\beta_{2,\alpha}\big(\eta_1-\beta_{1,\alpha})}{\eta_1\big(
\beta_{2,\alpha}-\beta_{1,\alpha}\big)}\, e^{-b \beta_{1,\alpha}}
+\frac{\beta_{1,\alpha}\big(\beta_{2,\alpha}-\eta_1\big)}
{\eta_1\big(\beta_{2,\alpha}-\beta_{1,\alpha}\big)}\,e^{-b \beta_{2,\alpha}}\!\!\bigg)
 e^{\alpha t}d\alpha,
\end{equation}
where $u>0$ is an arbitrary point such that the vertical line by $u$ does not pass through any  singular
 point mentioned in Proposition \ref{prop_lap_comp1}, and    $\beta_{1,\alpha}$ and $\beta_{2,\alpha}$
are the roots of $G(z)=\alpha=u+iv$ with
positive real part.
\end{proposition}

\begin{remark}
\label{remark:simetria}
From  a practical point of view 
 it is worth remarking that the expression
within the integral in (\ref{laplace:inversio}) is invariant by the permutation
of the roots $\beta_{1,\alpha}$ and $\beta_{2,\alpha}$.
 So, to compute numerically the Bromwich integral in (\ref{laplace:inversio})
  it suffices to
 compute the four roots of $P_\alpha(z)$ and take the  roots with positive
 real part.

\end{remark}
\section{Laplace transform of the joint distribution of the process and the running maxima}
Our objective is to compute the Laplace transform of the joint  probability
$$\p\{X_t\ge a,\, \tau_b\le t\}=\p(X_t\ge a, \max_{0\le s\le t} X_s\ge b)$$
 for $a\le b$ and $b>0$. The variable $\max_{0\le s\le t} X_s$ is called
 {\it the running maxima} of the process.
Our main result is the following:

\begin{proposition}
\label{prop:laplace:joint}
Fix $a\le b$ and $b>0$. Then for $\alpha\in \C$ with ${\rm Re}(\alpha)>0$ except
 at
the removable  singularities mentioned in Proposition \ref{prop_lap_comp1},
\begin{align}
\int_0^\infty   e^{-\alpha t}&\p\{X_t\ge a,\, \tau_b\le t\}\,dt
=\frac{1}{\alpha}\big(A(\alpha)+ B(\alpha)\big)
 \notag\\
& +\big(A(\alpha)\,C_3(\alpha)+B(\alpha)\,D_3(\alpha)\big) e^{-(b-a) \beta_{3,\alpha}}\notag\\
&+\big(A(\alpha)\,C_4(\alpha) +B(\alpha)\,D_4(\alpha)\big)e^{-(b-a)\beta_{4,\alpha}},
 \label{pro-laplace1}
\end{align}
where
\begin{equation}
\label{Aq}
A(\alpha)=\E\big[e^{-\alpha b}\bf 1_{\{X_{\tau_b}=b\}}\big]=\frac{\eta_1-\beta_{1,\alpha}}{\beta_{2,\alpha}-\beta_{1,\alpha}}\,
 e^{-b\beta_{1,\alpha} }
+\frac{\beta_{2,\alpha}-\eta_1}{\beta_{2,\alpha}-\beta_{1,\alpha}}\,
e^{- b\beta_{2,\alpha}},
\end{equation}
\begin{equation}
\label{Bq}
B(\alpha)=\E\big[e^{-\alpha b }\bf 1_{\{X_{\tau_b>b\}}}\big]=\frac{\big(\beta_{2,\alpha}-\eta_1\big)\big(\eta_1-\beta_{1,\alpha}\big)}
{\eta_1\big(\beta_{2,\alpha}-\beta_{1,\alpha}\big)}\,\Big(e^{-b\beta_{1,\alpha} }
-e^{-b\beta_{2,\alpha}}\Big),
\end{equation}
\begin{equation}
\label{C}
C_j(\alpha)=\frac{1}{\beta_{j,\alpha}G'(-\beta_{j,\alpha})},\ j=3,4,
\end{equation}
\begin{equation}
\label{D}
D_j(\alpha)=\frac{\eta_1}{(\eta_1+\beta_{j,\alpha})
\beta_{j,\alpha}G'(-\beta_{j,\alpha})},\ j=3,4,
\end{equation}
$G'(z)$ is the derivative of $G(z)$, and $-\beta_{3,\alpha}$
and $-\beta_{4,\alpha}$ are the  zeros of $G(x)=\alpha$
with negative real part.
\end{proposition}

\begin{remark}
As in Remark \ref{remark:simetria}, the expression in the right hand side of (\ref{pro-laplace1})
is invariant by the permutation between the roots with positive real
part, and between the roots with negative real part.
\end{remark}

\begin{remark}
\label{rees_Lap1}
Note that with these notations, equality (\ref{laplace:uni}) can be written as
$$\int_0^\infty e^{-\alpha t}\,\p(\tau_b\le t)\, dt=
\frac{1}{\alpha}\big(A(\alpha)+B(\alpha)\big).$$
\end{remark}

Also we prove  that the Laplace transform (\ref{pro-laplace1}) can be inverted by Bromwich integral:

\begin{proposition}
\label{prop:laplace:inversio2}
For any $a\le b$, $b>0$ and $t>0$,
\begin{align*}
\p\{X_t\ge a,\, \tau_b\le t\}=
\frac{1}{2\pi i}& \int_{u-i\infty}^{u+i\infty}
\bigg(\frac{1}{\alpha}\big(A(\alpha)+ B(\alpha)\big)\\
 &+\big(A(\alpha)\,C_3(\alpha)+B(\alpha)\,D_3(\alpha)\big)
 e^{-(b-a) \beta_{3,\alpha}}\\
&+\big(A(\alpha)\,C_4(\alpha) +B(\alpha)\,D_4(\alpha)\big)
e^{-(b-a)\beta_{4,\alpha}}
\bigg)
 e^{\alpha t}\, d\alpha,
\end{align*}
where $u>0$ is an arbitrary point such that the vertical line by $u$ does not pass through any  singular
 point mentioned in Proposition \ref{prop_lap_comp1},
$\beta_{1,\alpha}$ and $\beta_{2,\alpha}$ (respectively  $-\beta_{3,\alpha}$
 and $-\beta_{4,\alpha}$
are the roots of $G(z)=\alpha=u+iv$ with
positive real part (resp., negative part), where $G$ is given in (\ref{G}), and
$A(\alpha),\ B(\alpha),\,C_j(\alpha)$ and $D_j(\alpha)$, $j=3,4$ are given in
(\ref{Aq})-(\ref{D}).

\end{proposition}

\section{Numerical inversion of the Laplace transforms}
\subsection{Complex inversion}
\label{complexInversionNumerics}
We follow  the classical paper of Abate and Whitt (1995) (see also Abate and
 Whitt 1992).
In order to compute numerically a Bromwich integral of a Laplace transform $\wh f(\alpha)$,
$$f(t)=\frac{1}{2\pi i}\int_{u-i\infty}^{u+i\infty}\wh f(\alpha) e^{\alpha t}\,
 d\alpha,$$
  it is first rewritten as a real integral:
$$f(t)=\frac{2 e^{ut}}{\pi}\int_{0}^{\infty}{\rm Re}\big(\wh f(u+i x)\big) \cos(xt)\,
 dx.$$
This integral is numerically evaluated by means of the trapezoidal
rule. By taking a step size of length $\pi/(2t)$, half of the
 cosinus of the sum cancels and we get a nearly alternating series
$$f_d(t)=\frac{e^{A/2}}{2t}\, {\rm Re }(\wh f)\Big(\frac{A}{2t}\Big)+
\frac{e^{A/2}}{t}
\sum_{k=1}^\infty (-1)^k{\rm Re }(\wh f)\Big(\frac{A+2\pi  k i}{2t}\Big), $$
where $A=2t u.$
Using Poisson summation formula, the discretization error
can be controlled. Specifically, when the original function is bounded by 1, $\vert f(t)\vert\le 1$, as in the cases
that we are dealing, the discretization error can be bounded:
$$e_d:=\vert f(t)-f_d(t)\vert \le \frac{e^{-A}}{1-e^{-A}}\le e^{-A}.$$
Then,  we can take  $A=-\delta\,\log(10)$ to have at most $10^{-\delta}$
 discretization error.

Then the infinite sum is truncated to $n$ terms,
$$s_n(t)=\frac{e^{A/2}}{2t}\, {\rm Re }(\wh f)\Big(\frac{A}{2t}\Big)+
\frac{e^{A/2}}{t}
\sum_{k=1}^n {\rm Re }(\wh f)\Big(\frac{A+2\pi  k i}{2t}\Big),$$
and finally, as an  acceleration method, Euler summation is applied to $n$ terms
after an initial $B$:
$$E(n,B,t)=\sum_{k=0}^n \binom{n}{k}2^{-n}s_{B+k}(t).$$
To estimate the truncation error,  Abate and Whitt (1995)   suggest using
the difference of successive terms: $E(n,B,t)-E(n+1,B,t)$.

\subsection{Inversion of the real Laplace transform by Gaver-Stehfest method}
Kou and Wang (2003) invert the real Laplace transform (\ref{laplace:uni})
of Proposition \ref{Prop:teoKou} by using Gaver-Stehfest method.
This method is based on the Post-Widder formula
(see Abate and Whitt 1992, page 52) that says
that under some regularity conditions
$$f(t)=\lim_{n\to \infty} f_n(t),$$
where
$$f_n(t)=\frac{(-1)^n}{n!}\Big(\frac{n+1}{t}\Big)^{n+1}
\wh f^{(n)}\Big(\frac{n+1}{t}\Big),$$
where $\wh f^{(n)}$ is the $n$-derivative of the Laplace transform $\wh f$.
Gaver changes this derivative by a discrete analog by using finite differences;
specifically,
he proves  (Abate and Whitt 1992, Proposition 8.1)
that
$f(t)=\lim_{n\to \infty} \wt f_n(t),$
where
\begin{equation}
\label{factorials}
\wt f_n(t)=\frac{\log 2}{t}\frac{(2n)!}{n!(n-1)!}\sum_{k=0}^n(-1)^k\binom{n}{k}
\wh f\big((n+k)\log(2)/t\big).
\end{equation}
To accelerate the convergence, Kou and Wang (2003) use
   a variation of a method proposed by Stehfest (see Abate
  and Whitt 1992,  Proposition 8.2) computing
 $$f^*_n(t)=\sum_{k=1}^nw(k,n)\wt f_{B+k}(t),$$
 where
 $$w(k,n)=(-1)^{n-k}\,\frac{k^n}{k!(n-k)!},$$
 and $B$ is the initial burning-out number.

 As Abate and Whitt (1992)  page 52  remark, Gaver-Stehfest method is much
 less robust than  Fourier-series methods (such as the one based on Bromwich integral),
 and there is no error analysis. The computation of $f^*_n(t)$
is not easy due to the presence of the factorials in (\ref{factorials}).
 They recommend
taking $n=16$, and  using about 30 digits precision.

\section{A numeric example}
\label{sec_Num}
We analyze the example studied by Kou and Wang  (2003)  with
 parameters
$\mu=0.1, \,\sigma=0.2,\, p=0.5,\ \eta_1=50$, $\eta_2=100/3$ and $\lambda=3$.
Then
$G(z)$ given in (\ref{G}) is
$$G(z)=0.1z+0.02 z^2+3\Bigg(\frac{25}{50-z}+
 \frac{50/3}{100/3+z}-1\Bigg).$$
The zeros of $G(z)=\alpha$ are the roots of the polynomial (see \ref{GP})
$$P_\alpha(z)=-\frac{1}{50}\,{z}^{4}+ \frac{7}{30}\,{z}^{3}+
 (\alpha+38) {z}^{2}
+(-\frac{50}{3}\, \alpha+\frac{425}{3} ) z-\frac {5000}{3}\,\alpha,$$
The resultant of $P_\alpha$ and $\partial P_\alpha/\partial z$
is

\begin{align*}
R(\alpha)= -\frac{100}{9}\,&\alpha^5+\frac{244175}{162}\,\alpha^4-
\frac{115511483}{1458}\,\alpha^3
+\frac{3407354201}{1944}\,\alpha^2\\
&-\frac {17889068323}{972}\,\alpha-
\frac {8877228895}{5832},
\end{align*}
with zeros, approximately,
$$51.88+ 25.1\,i,\, 15.98+ 15.72\,i,\,- 0.08,\, 15.98- 15.72\,i,\, 51.88- 25.1\,i.$$
\begin{remark}The resultant of two polynomials can be computed easily in most computer
algebra systems, for example, in  Maple. Anyway,
a program in  C to compute the resultant can be asked to
Oscar Saleta (osr@mat.uab.cat)  from the Department of
Mathematics of the UAB.
\end{remark}

\bigskip

Like Kou and Wang (2003), we consider the values $t=1,\, b=0.3$ and
$a=0.2$. The times reported by Kou and Wang (2003) (with a
Pentium(R)
400 MHz PC) are given in Table \ref{kuo_wang_table}.
Given that the Laplace transform of the joint distribution provided
 by Kou and Wang
depends on the so-called {\sf Hh} functions, they need a large computing time with
MATHEMATICA.

\begin{table}[hbt]
\caption{Kou and Wang (2003) results. The accuracy in the computation of
$\Prob(\tau_b\le t)$ is not provided.}
\label{kuo_wang_table}
\centering
\begin{tabular}{ccc}
\toprule
& $\Prob(\tau_b\le t)$ & $\Prob(\tau_b\le t,\ X_t\ge a)$ \\
\midrule
Results &0.25584 & 0.22362\\
Accuracy & NP & $10^{-6}$\\
CPU time & 1.76 {\rm sec} & 4.61 {\rm min}\\
\bottomrule
\end{tabular}

\end{table}

To do the numeric inversion using Bromwich integral we use C language. With
the notations of
Subsection \ref{complexInversionNumerics}, we take $u=7,\ A=14,\ n=12$ and $B=4$.
The C programs were tested in an
Intel(R) Core(TM)2 Duo CPU E8400 @ 3.00GHz with 6GB of ram and Linux OS
 (Xubuntu 16.10 64bit).
The results are
given in the first two columns of  Table \ref{OurRes}. The accuracy is computed comparing
 the results given in table \ref{OurRes}
with the results for $n=50$.

\begin{table}[hbt]
\caption{Results using Bromwich integral with $n=12$ (first two columns) and Gaver-Stehfest
with 30 digits precision and  $n=10$ (third column) with a C program.}
\label{OurRes}
\centering
\begin{tabular}{cccc}
\toprule
& $\Prob(\tau_b\le t)$ & $\Prob(\tau_b\le t,\ X_t\ge a)$ & $\Prob(\tau_b\le t)$
  \\
& & &by Gaver-Stehfest\\
\midrule
Results &0.2558436 & 0.223616 & 0.2558433\\
Accuracy & $10^{-8}$ & $10^{-7}$\\
CPU time {\rm (ms)}&1.2  & 1.4  & 33.9 \\
\bottomrule
\end{tabular}

\end{table}

Unfortunately, Kou and Wang  (2003) do not provide information about
the program or language that they use to compute Gaver-Stehfest algorithm. With the purpose
of making a   clearer  comparison  we also
program Gaver-Stehfest algorithm with C  to estimate the probability $\Prob(\tau_b\le t)$
using 30 digits precision,
 and we get, with $n=10$, the result given in the third column of Table \ref{OurRes}.
We should point out that the precision of the Gaver-Stehfest algorithm
 is more difficult to control since if a  bigger
 $n$ is used  the algorithm starts to oscillate  and explodes.
The explosion time  can be delayed by increasing the number
 of digits in the computations and   in  Table \ref{Oscillation} there are
 some results using 30, 40 and 50 digits; this fact is a bit unsatisfactory.

  \begin{table}[hbt]
\caption{Results of Gaver-Stehfest algorithm with 30, 40 and 50
 digits precision when the number $n$ is increased.}
\label{Oscillation}
\centering
\begin{tabular}{ccccccc}
\toprule
Digits& \multicolumn{6}{c}{$n$}\\
\cmidrule(l){2-7}
 & 10 & 20 & 30 & 40 & 50 & 60\\
\cmidrule(l){2-7}
30 &  0.2558433 & 0.2558430 & 0.2558430 & 36238.016 &  -1.5949e17 & -1.0483e31\\
40  &   0.2558433 &  0.2558430 & 0.2558430 & 0.2558443 &  -3.5452e6 &  -2.2530e20\\
50 &   0.2558433 & 0.2558430 & 0.2558430 &  0.2558430 & 0.2537628 & -1.1717e10 \\
\bottomrule
\end{tabular}
\end{table}

Finally, as another benchmark, we also estimate by Monte Carlo
 simulation both probabilities
 $\Prob(\tau_b\le t)$ and $\Prob(\tau_b\le t,\ X_t\ge a)$. We should remark that
 the estimation of probabilities of first passage times of a jump process
 by Monte Carlo is a difficult
 problem that would require more research; see the comments by Kou and Wang
 (2003), pages 520 and 521. For present purposes we generate
  a discrete approximation of the Kou model along  a grid of 2000 points with
   {\sf \bfseries R} language, and we make $20000$ replications.
   The results, given in Table \ref{Montecarlo}, are very similar to the ones
   obtained by Kou and Wang
   (2003)  (last part of their Table 1), and, as they comment, the probabilities are underestimated.
   The CPU time (with  an Intel(R) Core2 Duo 3000 MHz PC)
   is 7.82 sec.
\begin{table}[hbt]
\caption{Computation of the probabilities by Monte Carlo simulation with a grid
of 2000 points and $20000$ replications.}
\label{Montecarlo}
\centering
\begin{tabular}{rccc}
\toprule
& $\Prob(\tau_b\le t)$ & $\Prob(\tau_b\le t,\ X_t\ge a)$\\
Point estimate & 0.25195 & 0.22105   \\
95\% IC & (0.2459,0.2580) & (0.2153,0.2268)\\
\bottomrule
\end{tabular}

\end{table}

\begin{remark}The codes of Maple programs to compute  the probabilities $\Prob(\tau_b\le t)$ and
 $\Prob(\tau_b\le t,\ X_t\ge a)$
 by Bromwich integral, and the former also by  Gaver-Stehfest algorithm, are given
  in  Appendix 2.  Also the C codes to compute both probabilities by
 Bromwich integral are in that Appendix. The get the C code to compute $\Prob(\tau_b\le t)$
 by  Gaver-Stehfest algorithm with multiprecision, the interested reader
 can write to Oscar Saleta (osr@mat.uab.cat).
\end{remark}
\section{Conclusion}
We have shown that the Laplace transform of the first passage times of
the Kou model can be  extended  to the complex half plane. We also give an explicit formula of
the complex Laplace transform of the joint distribution of the first passage time
and the running maxima without the use of the {\sf Hh} functions.
The  inversion of this last Laplace transform, compared with the method
given by Kou and Wang (2003),  produces a dramatic
reduction of the computing time.
In general  the use of complex Laplace transforms instead of  real
 Laplace transforms
has the advantage of improving the accuracy and provides a
 better control of the error.

\section*{Acknowledgements}
The authors are very grateful to  Oscar Saleta from the Department of
Mathematics of the UAB for the elaboration of the C programs; also
  would like to thank to Prof. Armengol Gasull from the Department of
Mathematics of the UAB for  fruitful conversations.
The second author was   partially
 supported   by grant
  MINECO reference  MTM2012-33937 and Generalitat de Catalunya
 reference 2014-SGR422.

\section*{Appendix 1: Proof of the propositions}

\subsection*{A1. Proof of Propositions \ref{prop_lap_comp1} and \ref{laplace:inv_complexa1}}
The following  property is well known; however,  our formulation
is not standard and  we believe it is convenient to write it  out.
For the sake of completeness  we include a
short proof.

\begin{lemma}
\label{lema:tranformada:p}
 Let $X$ be a non negative random variable. Then for every $\alpha\in \C$
 with $Re(\alpha)>0$,
\begin{equation}
\label{tranformada:p}
\E[e^{-\alpha X}]=\alpha\int_{0}^\infty e^{-\alpha x}\p\{X\le x\}\, dx.
\end{equation}
\end{lemma}

\noindent{\it Proof.} First we prove the equality for $\alpha>0$.
 The proof is based on the property that says that if $Z$ is
random variable such that $0\le Z\le 1$, then
$$\E[Z]=\int_0^1\p\{Z\ge t\}\, dt,$$
which follows from Fubini Theorem. From that property,
 for $\alpha>0$,
taking $Z=e^{-\alpha X}$ ,
 it is deduced (\ref{tranformada:p}).

 Now, since both sides of (\ref{tranformada:p})
  involve Laplace transforms well defined for $\alpha>0$, both parts
  can be extended to analytic functions on the half plane $\{\alpha\in \C: {\rm Re}(\alpha)>0\}$;
  since they coincide for $\alpha>0$, they coincide on the half plane.
   $\quad \blacksquare$

\bigskip

\noindent{\bf Proof of Proposition \ref{prop_lap_comp1}.}
As we commented after Proposition 1, $P_\alpha(z)$ is a four degree polynomial
whose coefficients are polynomials in $\alpha$. Then,  by analytic continuation,
we can  construct  four analytic   functions (in $\alpha$)
that give the zeros of this polynomial, see
 Ahlfors (1979), Chapter 8, Section 2.1.
 In principle, these analytic functions
exist for $\alpha$
  in any  part of the complex
  plane where
the  polynomial $P_\alpha(z)$ has  no multiple
zeros; since the multiple zeros of $P_\alpha(z)$ are zeros of both
$P_\alpha(z)$ and
$\partial P_\alpha/\partial z$, it suffices to look for the zeros of
the resultant $R(\alpha)$ of $P_\alpha(z)$ and $\partial P_\alpha/\partial z$,
which  is a polynomial in $\alpha$.
So the four functions are analytic except on  the zeros of $R(\alpha)$,
but at these points the singularities are removable (see the above mentioned
 reference to
 Ahlfors 1979).
 Hence, removing these singularities we get four analytic functions
 on the whole $\C$, that we call the global solutions;
  in this proof we denote these functions by
 $\bs \beta_{1}(\alpha),\, \bs \beta_{2}(\alpha),\,-\bs \beta_{3}(\alpha) $
 and $-\bs \beta_{4}(\alpha)$, where $\bs \beta_j(\alpha)=\beta_{j,\alpha}$
  when $\alpha>0$ (this happens by the continuity of the global solutions
  and the fact that for  $\alpha>0$ the four roots of $P_\alpha$
  are different, as we mentioned in Section \ref{sec_Kou}).

 Now, the expression on the right hand side of (\ref{laplace:uni}) can be extended
 to  $\C$ changing  $\beta_{j,\alpha}$  by $\bs\beta_j(\alpha)$; that is, define
\begin{equation}
\label{f1global}
 \wh f_1(\alpha)=
 \frac{1}{\alpha}\bigg(
 \frac{\bs\beta_2(\alpha)\big(\eta_1-\bs \beta_{1}(\alpha)\big)}{\eta_1\big(\bs \beta_{2}(\alpha)-
\bs\beta_{1}(\alpha)\big)}\, e^{-b \bs \beta_{1}(\alpha)}
+\frac{\bs\beta_1(\alpha)\big(\eta_1-\bs \beta_{2}(\alpha)\big)}
{\eta_1\big(\bs \beta_{1}(\alpha)-\bs \beta_{2}(\alpha)\big)}\, e^{-b \bs \beta_{2}(\alpha)}\bigg).
\end{equation}
This function is analytic except at $\alpha=0$ and for $\alpha$ such that
$\bs\beta_1(\alpha)=\bs \beta_2(\alpha).$ Now we prove that for these last points
the singularities are removable. Write $\wh f_1$
in the following way (we put $\beta_j$ instead of $\bs \beta_j(\alpha)$ to shorten the
notations):
$$\wh f_1=\frac{e^{\beta_1 b}}{\alpha \eta_1}\Big(\beta_1+
\frac{\beta_1(\beta_2-\eta_1)}{\delta}\big(e^{-\delta b}-1\big)\Big),$$
where $\delta=\beta_2-\beta_1$. Since $e^{-\delta b}-1=-\delta b+o(\delta),$
we have that $\lim_{\delta\to  0}f$ is finite, and so $f$ is bounded in a
 neighborhood
of $\delta=0$, which implies that the singularity is removable.

Moreover, for $\alpha\in\C$ with ${\rm Re}(\alpha)>0$,
the equation $G(z)=\alpha$ has no roots of the form $z=is$, for some real $s$;
this is due to the fact  that for any $s\in\R$, ${\rm Re}\big(G(is)\big)<0$,
as can be checked from the expression of $G(z)$ given in (\ref{G}).
This implies that for $\alpha\in{\C}$, with ${\rm Re}(\alpha)>0$  none of the functions
$\bs \beta_{1}(\alpha),\,\bs \beta_{2}(\alpha),\,-\bs\beta_{3}(\alpha)$ and
$-\bs \beta_{4}(\alpha)$
  crosses the imaginary axis.
As a consequence, we always have
 ${\rm Re}(\beta_{j}(\alpha))>0$, $j=1,\dots,4$.

On the other hand, since the Laplace transform on the left side of
 (\ref{laplace:uni}) is well defined for $\alpha>0$,
it can be extended to $\alpha\in \C$ with ${\rm Re}(\alpha)>0$, and it is an analytic
function on that half plane.
Since this extension and $\wh f_1$ coincide for $\alpha>0$
we get the desired result.   $\quad \blacksquare$

\bigskip

\noindent{\bf Proof of Proposition \ref{laplace:inv_complexa1}.}
A sufficient condition to apply the complex inversion formula to a
 function $\wh f$ analytic on $\C $, except on a finite number of
 isolated singularities all of them in a half plane
  $\{\alpha\in \C:\,{\rm Re} (\alpha) \le a\}$, is that there are positive constants
  $M, R, r$ such that for $\vert \alpha \vert >R,$
$$\vert\wh f (\alpha)\vert \le \frac{M}{\vert \alpha\vert^r},$$
see Marsden and Hoffman (1999) Corollary 8.2.2.

We want a bound of this type for the function $\wh f_1(\alpha)$  defined in (\ref{f1global}).
Extending  Kou and Wang (2003), Theorem 3.1, to the complex half plane
(or by  Lemma \ref{lema:tranformada:p} and Proposition \ref{Prop:teoKou})
it follows that for $\alpha\in\C$ with ${\rm Re}(\alpha)>0$,
$$\E[e^{-\alpha \tau_b}]=\alpha \int_0^\infty e^{-\alpha t}\,\p(\tau_b\le t)\, dt=
\alpha \wh f_1(\alpha),$$
 Hence,
$$\wh  f_1(\alpha)=\frac{1}{\alpha}\,\E[e^{-\alpha \tau_b}],$$
and since ${\rm Re}(\alpha)>0$ we have
$\big\vert \E[e^{-\alpha \tau_b}]\big\vert \le 1.$   $\quad \blacksquare$

\bigskip

\subsection*{A2. Proof of Propositions \ref{prop:laplace:joint} and \ref{prop:laplace:inversio2}}
We summarize in the following lemma the computation of the Fourier transform
of a rational function that we need:

\begin{lemma}
\label{lemma:fourier}
Let $R$ and $S$ be two real polynomials, with ${\rm deg}(R)<{\rm deg}(S)$, $R(0)\ne 0$,
and $S$ with only simple zeros $0,\gamma_1,\dots,\gamma_n$.
Then,
\smallskip

{\bf 1. } if $\omega<0$,
$$\int_{-\infty}^\infty e^{-i\omega s}\,\frac{R(is)}{S(is)}\,ds=
\pi \,\frac{R(0)}{S'(0)}+ 2\pi  \sum_{j: \,{\rm Re}(\gamma_j)<0}\frac{R(\gamma_j)}{S'(\gamma_j)} \,
 e^{-\omega \gamma_j}
$$

\smallskip

{\bf 2. } if $\omega>0$,
$$\int_{-\infty}^\infty e^{-i\omega s}\,\frac{R(is)}{S(is)}\,ds=
-\pi\,\frac{R(0)}{S'(0)}- 2\pi\sum_{j: \,{\rm Re}(\gamma_j)>0}\frac{R(\gamma_j)}{S'(\gamma_j)}\,
  e^{-\omega \gamma_j}
$$

\end{lemma}

\noindent{\it Proof.}

This is  just a particular case of the computation of the Fourier transform of
a rational function by calculus of residues. Define $\wt R(z)=R(iz)$ and
$\wt S(z)=S(iz)$. Then the rational fraction $\wt R/\wt S$ is analytic except
at the poles at 0, $-i\gamma_1,\dots, -i\gamma_n$. Applying standard formulas
(see, for example,   Marsden and Hoffman 1999,
Propositions 4.3.11 and 4.3.12)
to the computation of $\int_{-\infty}^\infty e^{-i\omega s}(\wt R(s)/\wt S(s))\,
  ds$ we get the result.$\quad \blacksquare$

\medskip

 A key point of the proof of Proposition \ref{prop:laplace:joint}
 is   the following proposition  from Kou and Wang (2003):

\begin{proposition} (Kou and Wang 2003, Proposition 4.1)
Fix $a\le b$ and $b>0$. For $\alpha>0$,
\begin{align}
\label{Prop-Kou}
\int_0^\infty & e^{-\alpha t}\p\{X_t\ge a,\, \tau_b\le t\}\,dt\notag\\
&=A(\alpha)\int_0^\infty  e^{-\alpha t}\p\{X_t\ge a-b\}\,dt+
B(\alpha)\int_0^\infty  e^{-\alpha t}\p\{X_t+\xi^+\ge a-b\}\,dt,
\end{align}
where $\xi^+$ is an exponential random variable with parameter $\eta_1$, independent
of the process $X$, and $A(\alpha)$ and $B(\alpha)$ are given in (\ref{Aq})
and (\ref{Bq}) respectively.

\end{proposition}

\noindent{\bf Proof of Proposition \ref{prop:laplace:joint}.}

\medskip

We first compute the Laplace transform for $\alpha>0$ and later we
extend the result to the complex half plane.

\medskip

\noindent{\bf A. Computation of the real Laplace transform}

\smallskip

\noindent{\bf 1.} To compute the first Laplace transform on the right
hand of (\ref{Prop-Kou}),
 consider,
for
$c<0<d$,
\begin{equation}
\label{laplace1}
\int_0^\infty e^{-\alpha t}\p\{c\le X_t\le d\}\,dt=
\int_{t=0}^\infty \int_{x=c}^de^{-\alpha t} f(t,x)\,dx\,dt,
\end{equation}
where $f(t,x)$ is the density function of $X_t$.
Denote by $\phi_X$ the characteristic function of a random variable $X$.
Since the Brownian motion $\{W_t,\, t\ge 0\}$ and the jumps part of $X_t$ are independent, the
 characteristic function of $X_t$ satisfies
$$\phi_{X_t}(s)=\E[e^{isX_t}]=e^{tG(is)}=\phi_{\sigma W_t}(s)\wt \phi(s),$$
where $\wt \phi$ is another characteristic function . Then
$$\big\vert \phi_{X_t}(s)\big\vert\le \big\vert \phi_{\sigma W_t}(s)\big\vert=
e^{-\sigma^2 s^2t/2},$$
which is integrable (remember that $\sigma>0$). Then, by the inversion formula of integrable characteristic
functions,
$$f(t,x)=\frac{1}{2\pi}\int_{s=-\infty}^\infty e^{-isx}e^{t G(is)}\,ds.$$
Hence,
$$(\ref{laplace1})=\frac{1}{2\pi}\int_{t=0}^\infty \int_{x=c}^d
\int_{s=-\infty}^\infty  e^{-\alpha t} e^{-isx}e^{tG(is)}\,ds
\,dx\,dt.$$
Furthermore,
\begin{align*}
\int_{t=0}^\infty\int_{x=c}^d\int_{s=-\infty}^\infty& \big\vert
e^{-\alpha t} e^{-isx}e^{tG(is)}\big\vert \,dt\,dx\,ds\\
&=(d-c)\int_{t=0}^\infty e^{-\alpha t}\Bigg( \int_{s=-\infty}^\infty
e^{-\sigma^2s^2t/2}\,ds\Bigg)dt\\
&=(d-c)\sigma^{-1}\sqrt{2\pi}\int_{t=0}^\infty \frac{1}{\sqrt t}\,e^{-\alpha t}\, dt=
(d-c)\sqrt 2 \pi \frac{1}{\sigma\sqrt \alpha}<\infty.
\end{align*}
Using Fubini Theorem we can choose the most convenient order to perform
the iterated integrals.
Therefore,
\begin{align}
\label{laplace_pas}
(\ref{laplace1})&=\frac{1}{2\pi}\int_{s=-\infty}^\infty\Bigg(\int_{t=0}^\infty
 e^{-\alpha t}
e^{tG(is)}\Big(
\int_{x=c}^d  e^{-isx}\,dx\Big)dt\Bigg)ds\notag\\
&=\frac{1}{2\pi }\int_{s=-\infty}^\infty
\frac{1}{is}\Big(e^{-ics}-e^{-ids}\Big)\Bigg(\int_{t=0}^\infty
e^{-t(\alpha-G(is))}\,dt\Bigg)ds\notag\\
&=\frac{1}{2\pi}\int_{s=-\infty}^\infty
\frac{1}{is\big(\alpha-G(is)\big)}\Big(e^{-ics}-e^{-ids}\Big)ds.
\end{align}
Now, use the notations  in (\ref{GP}):
$$\frac{1}{\alpha-G(is)}=-\frac{Q(is)}{P_\alpha(is)},$$
and apply Lemma \ref{lemma:fourier} with
$$R(z)=-Q(z)\quad {\rm and}\quad S(z)=zP_\alpha(z).$$
Note that from (\ref{GP}),
$$G'(z)=\frac{P'_\alpha(z)Q(z)-P_\alpha(z)Q'(z)}{Q^2(z)},$$
and hence,
$$G'(\beta_j)=\frac{P'_\alpha(\beta_j)}{Q(\beta_j)}, \, j=1,2,$$
and
\begin{equation}
\label{Gprimaaval}
G'(-\beta_j)=\frac{P'_\alpha(-\beta_j)}{Q(-\beta_j)}, \, j=3,4.
\end{equation}
Combining all these formulas, we deduce that
\begin{align}
\label{laplace_pas2}
(\ref{laplace_pas})=\frac{1}{\alpha}+
C_3 e^{c \beta_{3}}+
C_4 e^{c\beta_{4}}+C_1e^{-d \beta_{1}}+C_2 e^{-d \beta_{2}},
\end{align}
where $$C_j=C_j(\alpha)=\frac{1}{\beta_{j,\alpha}G'(-\beta_{j,\alpha})},\ j=3,4,$$
 and  $C_1$ and $C_2$ do not matter.

\medskip

\noindent{\bf 2.} Now we have that
\begin{align*}
\int_0^\infty  e^{-\alpha t}\p\{X_t\ge a-b\}\,dt&=
\int_0^\infty  e^{-\alpha t}\Big(\lim_{d\to\infty}\p\{d\ge X_t\ge a-b\}\Big)\,dt\\
&=\lim_{d\to\infty}\int_0^\infty  e^{-\alpha t}\p\{d\ge X_t\ge a-b\}\,dt,
\end{align*}
where the last equality follows from  dominated convergence Theorem.
Passing to the limit (\ref{laplace_pas2}), and taking
 $c:=a-b<0$, we arrive at expression (\ref{pro-laplace1}).

 \bigskip

 The second Laplace transform in (\ref{Prop-Kou}) is computed exactly in the
 same way, using that the characteristic function of $X_t+\xi^+$ is
 $$\phi_{X_t+\xi^+}(s)=\phi_{X_t}(s)\, \frac{\eta_1}{\eta_1-i s}.$$

\medskip

\noindent{\bf B. Extension to the Laplace transform to the complex half plane.}  In the same way that in the
proof of Proposition \ref{prop_lap_comp1}, the right hand side of
(\ref{pro-laplace1}) can be extended to the complex plane except at some
isolated singularities. Denote this function by $\wh f_2(\alpha)$. We have
that
$$\wh f_2(\alpha)= \wh f_1(\alpha)+e^{-c\beta_3}\Big[\frac{F_3(\alpha)}{\beta_3 G'(-\beta_3)}
+\frac{F_4(\alpha)}{\beta_4 G'(-\beta_4)}\, e^{-c(\beta_4-\beta_3)}\Big],$$
where $\wh f_1(\alpha)$ is the function introduced in the proof of
Proposition \ref{prop_lap_comp1} (see Remark \ref{rees_Lap1}), $c=b-a>0$ and
$$F_j(\alpha)=A(\alpha)+B(\alpha)\,\frac{\eta_1}{\eta_1+\beta_j}, \ j=3,4.$$
Given that $\eta_1>0$ and ${\rm Re}(\beta_j)>0$, the function $\wh f_2$ only has singularities
where $G'(-\beta_j)=0,\ j=3,4$. However, by (\ref{Gprimaaval}),
$$G'(-\beta_j)=\frac{P'_\alpha(-\beta_j)}{Q(-\beta_j)}, \, j=3,4,$$
and hence $-\beta_j$ should be a root of both $P_\alpha$ and $P_\alpha'$, and then
$\beta_3=\beta_4$. As in the proof of Proposition \ref{prop_lap_comp1}, it can be
proved that these singularities are removable. $\qquad \blacksquare$

\bigskip

\noindent{\bf Proof of Proposition \ref{prop:laplace:inversio2}.}

\medskip

We are going to provide a bound of the type given in the
 proof of Proposition \ref{laplace:inv_complexa1} for the function $\wh f_2(z)$
  defined
on  the right hand side of (\ref{pro-laplace1}):
\begin{align*}
 \label{pro-laplace2}
\wh f_2(z)&=\frac{1}{\alpha}\big(A(\alpha)+ B(\alpha)\big)
 +\big(A(\alpha)\,C_3(\alpha)+B(\alpha)\,D_3(\alpha)\big)
  e^{-(b-a) \beta_{3,\alpha}} \\
&\qquad +\big(A(\alpha)\,C_4(\alpha) +B(\alpha)\,D_4(\alpha)\big)e^{-(b-a)\beta_{4,\alpha}}.
\end{align*}
Since
$A(\alpha)=\E\big[e^{-\alpha b}\bf 1_{\{X_{\tau_b}=b\}}\big]$ and
$B(\alpha)=\E\big[e^{-\alpha b}\bf 1_{\{X_{\tau_b>b\}}}\big]$ are bounded by 1,   we need only to find
 bounds for $C_j$ and $D_j$, $j=3,4$. These bounds are deduced from the following estimations
of the roots $\beta_3$ and $\beta_4$:
There is $R>0$  such that for $\vert \alpha\vert >R$, $\vert \beta_4\vert > C \alpha ^{1/4}$
and $\vert -\beta_3+\eta_2\vert< C /\vert \alpha\vert$. These estimations are
 proved below.

  From these estimations,
it is clear that for $\vert \alpha\vert >R$, there is $C>0$ such that
$$ \vert C_4(\alpha)\vert=
\Big\vert\frac{1}{\beta_{4}G'(-\beta_{4})}\Big \vert\le C \frac{1}{\vert \alpha \vert^{1/4}}$$
and
$$\vert D_4(\alpha)\vert=\Big\vert\frac{\eta_1}{(\eta_1+\beta_{4})
\beta_{4}G'(-\beta_{4})}\Big\vert \le C \frac{1}{\vert \alpha \vert^{1/4}}.$$

In relation to $C_3(\alpha)$ and $D_3(\alpha)$, note that in the expression
(\ref{GP}), $Q(z)=(\eta_1-z)(\eta_2+z)$, and, see (\ref{Gprimaaval}),
\begin{equation*}
G'(-\beta_3)=\frac{P'_\alpha(-\beta_3)}{Q(-\beta_3)}.
\end{equation*}
Hence, there is $R>0$ and $C>0$ such that for $\vert \alpha\vert>R$,
$$\vert C_3(\alpha)\vert<C/\vert \alpha \vert
\quad {\rm and} \quad \vert D_3(\alpha)\vert<C/\vert \alpha \vert.
$$

\bigskip

 The estimations of the roots $\beta_3$ and $\beta_4$ are deduced from Rouch\'e's Theorem,
 which says that if $f$ and $g$ are two analytic functions,
  $\gamma$  is a closed curve homotopic to a point, and on $\gamma$
 $$\vert f(z)-g(z)\vert  < \vert f(z) \vert,$$
 then $f$ and $g$ have the same number of zeros enclosed by $\gamma$
 (see, for example, Ahlfors 1999, Pag. 153).

To estimate $\beta_3$, note  that the polynomial $P_\alpha(z)$ has the form
 $$P_\alpha(z)=a_4 z^4+a_3 z^3+(a_2\alpha+a'_2 )z^2+(a_1\alpha+a'_1 )z+
a_0\alpha,$$
where all $a_j$ and $a'_j$ are constants. Then, applying Rouché's Theorem to
$f(z)=(a_2\alpha+a'_2 )z^2$,  $g(z)=f(z)-P_\alpha(z)$, and $\gamma$ the circle
centered at 0 of radius $C \alpha ^{1/4}$, where $C$ is a constant depending
on  $a_j$ and $a'_j$, for $\vert \alpha\vert$ big enough, we deduce that $P_\alpha(z)$ has two roots inside that
circle and two outside. Hence $\vert \beta_4\vert > C \alpha ^{1/4}$.

To estimate $\beta_3$, apply Rouch\'e's Theorem to
 $f(z)=\lambda(1-p)\eta_2/(\eta_2+z)-\alpha$,
$g(z)=f(z)-G(z)$ and $\gamma$ a circle centered at $-\eta_2$ (remember that it is
a real number) and radius $C/\vert \alpha\vert$ for $\vert \alpha\vert$ big enough. $\blacksquare$

\section*{Appendix 2. Maple and C codes}

\section*{Maple Codes}

\subsection*{Code 1: Inversion of the Laplace transform of $\bs{\Prob\{\tau_b\le t\}}$
by Bromwich integral: function { \tt\bfseries f1}}

\begin{Verbatim}[frame=lines,numbers=left,fontsize=\footnotesize]
Parameters:= proc(mu_, sigma_, lambda_, eta1_, eta2_, p_)
  global mu, sigma, lambda, eta1, eta2, p;
  mu := mu_;
  sigma := sigma_;
  lambda := lambda_;eta1 := eta1_;
  eta2 := eta2_;
  p := p_;
end proc:

Parameters(1/10,2/10,3,50,100/3,1/2):

Polynomial:= proc()
  global P;
  P := (-lambda+mu*z+(1/2)*sigma^2*z^2)*(eta1-z)*(eta2+z)+lambda*(p*eta1*(eta2+z)
          +(1-p)*eta2*(eta1-z));
  P := -2*collect(P-(eta1-z)*(eta2+z)*alpha, z)/sigma^2;
  P := unapply(P, z, alpha);
  return;
end proc:

Polynomial();

suma:=proc(f,t,A,n)
  local sum_a_k,s_n;
  sum_a_k:=add(evalf((-1)^k*Re(f((A+2*k*Pi*I)/(2*t)))),k=1..n);
  s_n:=exp(A/2)/(2*t)*Re(f(A/(2*t)))+exp(A/2)/t*sum_a_k;
  return evalf(s_n);
end proc:

Euler:=proc(f,t,A,n,B)
  local E,k;
  E:=2^(-n)*suma(f,t,A,B);
  for k from 1 to n do
  E:=E+2^(-n)*binomial(n,k)*suma(f,t,A,B+k);
  end do;
end proc:

f1:=proc(t,b,A,n,B)
  local hat_f1;
  hat_f1:=proc(alpha)
    local beta,beta1,beta2,denom,num; Polynomial():
    beta:=fsolve(P(z,alpha)=0,z);
    beta1:=beta[3];
    beta2:=beta[4];
    denom:=alpha*eta1*(beta2-beta1);
    num:=beta2*(eta1-beta1)*exp(-beta1*b)+beta1*(beta2-eta1)*exp(-beta2*b);
    return evalf(num/denom);
  end proc;
  Euler(hat_f1,t,A,n,B);
end proc:

f1(1,0.3,14,12,4);
\end{Verbatim}

Lines 1 to 8  fix the parameters of the model. Lines 12 to 19
 compute
a polynomial with the same roots as $P_\alpha$. In lines 23 to 28 there is a procedure
to compute   the discretization
of an integral. In lines 30 to 36 the acceleration by Euler summation
 is computed.
Finally, lines 38 to 50 introduce the Laplace transform to be inverted and compute the
inversion.

\subsection*{Code 2: Inversion of the Laplace transform of ${\bs
\Prob\{\ X_t\ge a, \,\tau_b\le t\}}$
by Bromwich integral: function {\bfseries\tt f2}}

\begin{Verbatim}[frame=lines,numbers=left,fontsize=\footnotesize]

Lines 1 to 36 of Code 1

f2:=proc(t,a,b,A,n,B)
  local hat_f2;
  hat_f2:=proc(alpha)
    local G,dG,beta,beta1,beta2,beta3,beta4,aux,Aalpha,Balpha,C3,C4,D3,D4,c;
    G:=z->mu*z+1/2*sigma^2*z^2+lambda*(p*eta1/(eta1-z)+(1-p)*eta2/(eta2+z)-1);
    dG:=D(G);
    beta:=fsolve(P(z,alpha)=0,z);
    beta3:=-beta[2];
    beta4:=-beta[1];
    beta1:=beta[3];
    beta2:=beta[4];
    aux:=(beta1,beta2)->(eta1-beta1)*exp(-beta1*b)/(beta2-beta1);
    Aalpha:=aux(beta1,beta2)+aux(beta2,beta1);
    Balpha:=aux(beta1,beta2)*(beta2-eta1)/eta1-aux(beta2,beta1)*(eta1-beta1)/eta1;
    C3:=1/(beta3*dG(-beta3));
    C4:=1/(beta4*dG(-beta4));
    D3:=eta1/((eta1+beta3)*beta3*dG(-beta3));
    D4:=eta1/((eta1+beta4)*beta4*dG(-beta4));
    c:=a-b;
    return (Aalpha+Balpha)/alpha+(C3*Aalpha+D3*Balpha)*exp(c*beta3)
               +(C4*Aalpha+D4*Balpha)*exp(c*beta4);
  end proc;
  Euler(hat_f2,t,A,n,B);
end proc:

f2(1,0.2,0.3,14,12,4);
\end{Verbatim}
\subsection*{Code 3:  Inversion of the Laplace transform of
${\bs\Prob\{\tau_b\le t\}}$
by  Gaver-Stehfest method: function {\tt f3}}

\begin{Verbatim}[frame=lines,numbers=left,fontsize=\footnotesize]
Digits:=30:

Lines 1 to 21 of the code 1

fn:=proc(f,n,t)
  local suma, k;
  suma:=0;
  for k from 0 to n do
  suma:=suma+ln(2)*(2*n)!/((n-1)!*k!*(n-k)!*t)*(-1)^k*f((n+k)*ln(2)/t);
  end do;
  evalf(suma);
end proc:

gaver:=proc(f,t,n,B)
  local suma, k;
  suma:=0;
  for k from 1 to n do
  suma:=suma+(-1)^(n-k)*k^n/(k!*(n-k)!)*fn(f,k+B,t);
  end do;
  evalf(suma);
end proc:


f3:=proc(t,b,n,B)
  local hat_f1;
  hat_f1:=proc(alpha)
    local beta,beta1,beta2,denom,num;
    beta:=fsolve(P(z,alpha)=0,z);
    beta1:=beta[3];
    beta2:=beta[4];
    denom:=alpha*eta1*(beta2-beta1);
    num:=beta2*(eta1-beta1)*exp(-beta1*b)+beta1*(beta2-eta1)*exp(-beta2*b);
    return evalf(num/denom);
  end proc;
  gaver(hat_f1,t,n,B);
end proc:

f3(1,0.3,10,2);
\end{Verbatim}

\section*{C Codes}

\subsection*{Code 4: Inversion of the Laplace transform of $\bs{\Prob\{\tau_b\le t\}}$
by Bromwich integral: function {\tt f1}}

\begin{Verbatim}[frame=lines, numbers=left,fontsize=\footnotesize]
#include "koujdm_lib.h"

#include <stdio.h>
#include <stdlib.h>
#include <inttypes.h>

int main(int argc, char *argv[]) {
    Parameters prm;
    Parametersf1 *prm2 = (Parametersf1 *)malloc(sizeof(Parametersf1));

    if (argc != 12
        || sscanf(argv[1],"%Lf",&prm.mu)!=1
        || sscanf(argv[2],"%Lf",&prm.sigma)!=1
        || sscanf(argv[3],"%Lf",&prm.lambda)!=1
        || sscanf(argv[4],"%Lf",&prm.eta1)!=1
        || sscanf(argv[5],"%Lf",&prm.eta2)!=1
        || sscanf(argv[6],"%Lf",&prm.p)!=1
        || sscanf(argv[7],"%"PRIi64,&prm2->t)!=1
        || sscanf(argv[8],"%Lf",&prm2->b)!=1
        || sscanf(argv[9],"%"PRIi64,&prm2->A)!=1
        || sscanf(argv[10],"%"PRIi64,&prm2->n)!=1
        || sscanf(argv[11],"%"PRIi64,&prm2->B)!=1) {
        fprintf(stderr,"%s mu sigma lambda eta1 eta2 p t b A n B\n",argv[0]);
        return -1;
    }

    long double result = euler(&hat_f1,prm2->t,prm2->A,prm2->n,prm2->B,prm,prm2);
    printf("result=%30.29Lg\n",result);

    return 0;
}
\end{Verbatim}

\subsection*{Code 5: Inversion of the Laplace transform of ${\bs
\Prob\{\ X_t\ge a, \,\tau_b\le t\}}$
by Bromwich integral. Function {\tt f2}}

\begin{Verbatim}[frame=lines,numbers=left,fontsize=\footnotesize]
#include "koujdm_lib.h"

#include <stdio.h>
#include <stdlib.h>
#include <inttypes.h>

int main(int argc, char *argv[]) {
    Parameters prm;
    Parametersf2 *prm2 = (Parametersf2 *)malloc(sizeof(Parametersf2));

    if (argc != 13
        || sscanf(argv[1],"%Lf",&prm.mu)!=1
        || sscanf(argv[2],"%Lf",&prm.sigma)!=1
        || sscanf(argv[3],"%Lf",&prm.lambda)!=1
        || sscanf(argv[4],"%Lf",&prm.eta1)!=1
        || sscanf(argv[5],"%Lf",&prm.eta2)!=1
        || sscanf(argv[6],"%Lf",&prm.p)!=1
        || sscanf(argv[7],"%"PRIi64,&prm2->t)!=1
        || sscanf(argv[8],"%Lf",&prm2->a)!=1
        || sscanf(argv[9],"%Lf",&prm2->b)!=1
        || sscanf(argv[10],"%"PRIi64,&prm2->A)!=1
        || sscanf(argv[11],"%"PRIi64,&prm2->n)!=1
        || sscanf(argv[12],"%"PRIi64,&prm2->B)!=1) {
        fprintf(stderr,"%s mu sigma lambda eta1 eta2 p t a b A n B\n",argv[0]);
        return -1;
    }

    long double result = euler(&hat_f2,prm2->t,prm2->A,prm2->n,prm2->B,prm,prm2);
    printf("result=%30.29Lg\n",result);

    return 0;
}
\end{Verbatim}

\subsection*{Code 6: Library {\tt koujdm\_lib.c}}

\begin{Verbatim}[frame=lines,numbers=left,fontsize=\footnotesize]
#include "koujdm_lib.h"

#include <inttypes.h>

/****
 * GENERIC FUNCTIONS
 ****/
// binomial coefficient computed recursively
int64_t binomial(int64_t n, int64_t k) {
    if (n<k)
        return -1;
    if (k==0)
        return 1;
    if (k > n/2)
        return binomial(n,n-k);
    return n*binomial(n-1,k-1)/k;
}

// coefs. of the polynomial of degree 4
Quartic polynomial(long double complex alpha, Parameters prm) {
    Quartic p;
    long double sigma2 = prm.sigma*prm.sigma;
    p.a = 1;
    p.b = -prm.eta1 + prm.eta2 + 2.*prm.mu/sigma2;
    p.c = -prm.eta1*prm.eta2 + 2./sigma2*(-prm.mu*prm.eta1 - prm.lambda + prm.mu*prm.eta2 - alpha);
    p.d = 2./sigma2*(-prm.eta1*prm.eta2*prm.mu - prm.eta1*prm.lambda*(prm.p-1)
        - prm.eta2*prm.lambda*prm.p + alpha*(prm.eta1-prm.eta2));
    p.e = 2.*prm.eta1*prm.eta2*alpha/sigma2;

    return p;
}

// direct computation of the roots of the polynomial of degree 4
//http://math.stackexchange.com/a/786
QuarticSolutions quartic_solve(Quartic p) {
    long double complex p1 = 2*p.c*p.c*p.c-9*p.b*p.c*p.d+27*p.a*p.d*p.d+27*p.b*p.b*p.e-72*p.a*p.c*p.e;
    long double complex p2 = p1+csqrt(-4*cpow((p.c*p.c-3*p.b*p.d+12*p.a*p.e),3.)+p1*p1);
    long double complex p3 = (p.c*p.c-3*p.b*p.d+12*p.a*p.e)/(3*p.a*cpow(p2/2.,1/3.))
        +cpow(p2/2.,1/3.)/(3*p.a);
    long double complex p4 = csqrt(p.b*p.b/(4*p.a*p.a)-2*p.c/(3*p.a)+p3);
    long double complex p5 = p.b*p.b/(2*p.a*p.a)-4*p.c/(3*p.a)-p3;
    long double complex p6 = (-cpow(p.b/p.a,3.)+4*p.b*p.c/(p.a*p.a)-8*p.d/p.a)/(4*p4);
    QuarticSolutions result;
    result.x1 = -p.b/(4*p.a)-p4/2-csqrt(p5-p6)/2;
    result.x2 = -p.b/(4*p.a)-p4/2+csqrt(p5-p6)/2;
    result.x3 = -p.b/(4*p.a)+p4/2-csqrt(p5+p6)/2;
    result.x4 = -p.b/(4*p.a)+p4/2+csqrt(p5+p6)/2;
    return result;
}

/*
 * FIRST BLOCK (f1, f2)
 */
// used in f1.c and f2.c
long double suma(long double complex (*f)(long double complex, Parameters, void*), int64_t t,
       int64_t A, int64_t n, Parameters prm, void* prm2) {
    long double sum_a_k = 0;
    int64_t sign;
    for (int_fast64_t k=1;k<=n;k++) {
        sign = k%2 ? -1 : 1;
        sum_a_k += sign*creal(f((A+2.*k*M_PI*I)/(2.*t),prm,prm2));
    }
    long double s_n = exp(A/2.)/(2.*t)*creal(f(A/(2.*t),prm,prm2))+exp(A/2.)/t*sum_a_k;
    return s_n;
}

// used in f1.c and f2.c
long double euler(long double complex (*f)(long double complex, Parameters, void*), int64_t t,
      int64_t A, int64_t n, int64_t B, Parameters prm, void* prm2) {
    long double E = pow(2.,(long double)(-n))*suma(f,t,A,B,prm,prm2);
    for (int_fast64_t k=1;k<=n;k++) {
        E += pow(2.,(long double)(-n))*binomial(n,k)*suma(f,t,A,B+k,prm,prm2);
    }
    return E;
}

// used in f1.c
long double complex hat_f1(long double complex alpha, Parameters prm, void *prm2) {
    Parametersf1 *prm2_ = (Parametersf1*)(prm2);

    Quartic p = polynomial(alpha, prm);
    QuarticSolutions beta = quartic_solve(p);

    long double complex max1,min1,max2,min2;

    if (creal(beta.x1)>creal(beta.x2)) {
        max1 = beta.x1;
        min1 = beta.x2;
    } else {
        max1 = beta.x2;
        min1 = beta.x1;
    }
    if (creal(beta.x3)>creal(beta.x4)) {
        max2 = beta.x3;
        min2 = beta.x4;
    } else {
        max2 = beta.x4;
        min2 = beta.x3;
    }

    long double complex beta1,beta2;
    if (creal(max1)<creal(min2)) {
        beta1 = min2;
        beta2 = max2;
    } else if (creal(max2)<creal(min1)) {
        beta1 = min1;
        beta2 = max1;
    } else if (creal(max1)<creal(max2)) {
        beta1 = max1;
        beta2 = max2;
    } else {
        beta1 = max2;
        beta2 = max1;
    }

    long double complex denom = alpha*prm.eta1*(beta2-beta1);
    long double complex num = beta2*(prm.eta1-beta1)*cexp(-beta1*prm2_->b)
           +beta1*(beta2-prm.eta1)*cexp(-beta2*prm2_->b);
    return num/denom;
}

// used in f2.c
long double complex aux(long double complex beta1, long double complex beta2, long double eta1,
       long double b) {
    return (eta1-beta1)*cexp(-beta1*b)/(beta2-beta1);
}

// used in f2.c
long double complex dG(long double complex z, Parameters prm) {
    return prm.mu+prm.sigma*prm.sigma*z+prm.lambda*(prm.p*prm.eta1/((prm.eta1-z)*(prm.eta1-z))
         -(1-prm.p)*prm.eta2/((prm.eta2+z)*(prm.eta2+z)));
}

// used in f2.c
long double complex hat_f2(long double complex alpha, Parameters prm, void *prm2) {
    Parametersf2 *prm2_ = (Parametersf2*)(prm2);

    Quartic p = polynomial(alpha, prm);
    QuarticSolutions beta = quartic_solve(p);

    long double complex max1,min1,max2,min2;

    if (creal(beta.x1)>creal(beta.x2)) {
        max1 = beta.x1;
        min1 = beta.x2;
    } else {
        max1 = beta.x2;
        min1 = beta.x1;
    }
    if (creal(beta.x3)>creal(beta.x4)) {
        max2 = beta.x3;
        min2 = beta.x4;
    } else {
        max2 = beta.x4;
        min2 = beta.x3;
    }

    long double complex beta1,beta2,beta3,beta4;
    if (creal(max1)<creal(min2)) {
        beta1 = min2;
        beta2 = max2;
        beta3 = -max1;
        beta4 = -min1;
    } else if (creal(max2)<creal(min1)) {
        beta1 = min1;
        beta2 = max1;
        beta3 = -max2;
        beta4 = -min2;
    } else if (creal(max1)<creal(max2)) {
        beta1 = max1;
        beta2 = max2;
        if (creal(-min1)<creal(-min2)) {
            beta3 = -min2;
            beta4 = -min1;
        } else {
            beta3 = -min1;
            beta4 = -min2;
        }
    } else {
        beta1 = max2;
        beta2 = max1;
        if (creal(-min1)<creal(-min2)) {
            beta3 = -min2;
            beta4 = -min1;
        } else {
            beta3 = -min1;
            beta4 = -min2;
        }
    }

    complex long double Aalpha, Balpha;
    Aalpha = aux(beta1,beta2,prm.eta1,prm2_->b)+aux(beta2,beta1,prm.eta1,prm2_->b);
    Balpha = aux(beta1,beta2,prm.eta1,prm2_->b)*(beta2-prm.eta1)/prm.eta1
        -aux(beta2,beta1,prm.eta1,prm2_->b)*(prm.eta1-beta1)/prm.eta1;

    long double complex C3,C4,D3,D4;
    C3 = 1./(beta3*dG(-beta3,prm));
    C4 = 1./(beta4*dG(-beta4,prm));
    D3 = prm.eta1/((prm.eta1+beta3)*beta3*dG(-beta3,prm));
    D4 = prm.eta2/((prm.eta2+beta4)*beta4*dG(-beta4,prm));

    long double c = prm2_->a-prm2_->b;

    return (Aalpha+Balpha)/alpha + (C3*Aalpha+D3*Balpha)*cexp(c*beta3)
             +(C4*Aalpha+D4*Balpha)*cexp(c*beta4);
}
\end{Verbatim}

\subsection*{Code 7: Library {\tt koujdm\_lib.h}}

\begin{Verbatim}[frame=lines,numbers=left,fontsize=\footnotesize]
#pragma once

#include <stdio.h>
#include <complex.h>
#include <math.h>
#include <stdarg.h>
#include <stdint.h>
#include <gmp.h>
#include <mpfr.h>
#include <mpc.h>

#define MPCRND MPC_RNDNN
#define MPFRND MPFR_RNDN


#ifndef M_PI
#define M_PI acos(-1.0)
#endif

/* Structures */

// quartic polynomial
typedef struct {
    long double complex a,b,c,d,e;
} Quartic;

typedef struct {
    long double complex x1,x2,x3,x4;
} QuarticSolutions;



// generic parameter structure
typedef struct {
    long double mu,sigma,lambda,eta1,eta2,p;
} Parameters;

// specific parameter structures
typedef struct {
    long double b;
    int64_t t,A,n,B;
} Parametersf1;

typedef struct {
    long double a,b;
    int64_t t,A,n,B;
} Parametersf2;




/* Generic functions */
// binomial
int64_t binomial(int64_t, int64_t);
// quartic polynomial
Quartic polynomial(long double complex, Parameters);
// roots of quartic polynomial
QuarticSolutions quartic_solve(Quartic);


// first block functions
long double suma(long double complex (*f)(long double complex, Parameters, void*), int64_t, int64_t,
       int64_t, Parameters, void*);
long double euler(long double complex (*f)(long double complex, Parameters, void*), int64_t, int64_t,
       int64_t, int64_t, Parameters, void*);
// f1.c
long double complex hat_f1(long double complex, Parameters, void*);
// f2.c
long double complex aux(long double complex, long double complex, long double, long double);
long double complex dG(long double complex, Parameters);
long double complex hat_f2(long double complex, Parameters, void*);
\end{Verbatim}

\newpage

\subsection*{Code 8: Makefile}

\begin{Verbatim}[frame=lines,numbers=left,fontsize=\footnotesize]
CFLAGS= -std=c11 -Wall -O3
LFLAGS= -lm -lgmp -lmpfr -lmpc

target: all

all: f1 f2

# Main executable
f1: f1.o koujdm_lib.o
    gcc $(CFLAGS) f1.o koujdm_lib.o -o f1 $(LFLAGS)

f2: f2.o koujdm_lib.o
    gcc $(CFLAGS) f2.o koujdm_lib.o -o f2 $(LFLAGS)



# Objects compilation
f1.o: f1.c
    gcc $(CFLAGS) -c f1.c $(LFLAGS)

f2.o: f2.c
    gcc $(CFLAGS) -c f2.c $(LFLAGS)


# Run configuration
runf1: f1
    ./f1 0.1 0.2 3 50 33.3333333333333333 0.5 1 0.3 14 12 4

runf2: f2
    ./f2 0.1 0.2 3 50 33.3333333333333333 0.5 1 0.2 0.3 14 12 4


# Cleaning directives
clean:
rm *o f1 f2
\end{Verbatim}

\end{document}